\documentclass[reqno,11pt,a4paper]{amsart}
\usepackage[utf8]{inputenc}

\usepackage{amsmath}
\usepackage{amsthm}
\usepackage{amssymb}
\usepackage{hyperref}
\usepackage{theoremref}
\usepackage{bbm}
\usepackage{etoolbox}
\usepackage{mathtools}
\usepackage{comment}
\usepackage[foot]{amsaddr}
\usepackage{enumitem}

\usepackage{tikz}
\usetikzlibrary{matrix,arrows,calc,cd,decorations.markings}

\input diagxy

\usepackage{scalerel,lmodern}
\input pdf-trans
\newbox\qbox
\def\usecolor#1{\csname\string\color@#1\endcsname\space}
\newcommand\outline[1]{\leavevmode%
  \def\maltext{#1}%
  \setbox\qbox=\hbox{\maltext}%
  \boxgs{Q q 2 Tr \thickness\space w 0 0 0 rg 0 G}{}%
  \copy\qbox%
}

\makeatletter
\renewcommand{\do}[1]{\@namedef{c#1}{\ensuremath{\mathcal{#1}}}}
\makeatother
\docsvlist{A,B,C,D,E,P,U,V} 

\newcommand{\lbar}{\underline{\hspace{0.8em}}}

\DeclareMathOperator{\id}{id}

\renewcommand*{\&}{%
  \relax
  \ifmmode
    \mathbin{\char`\&}%
  \else
    \char`\&\relax
  \fi
}

\renewcommand{\Im}{\mathrm{Im}}

\newtheorem{theorem}{Theorem}[section]
\theoremstyle{plain}

\newtheorem{corollary}[theorem]{Corollary}
\newtheorem{lemma}[theorem]{Lemma}

\numberwithin{equation}{section}
\theoremstyle{definition}

\newtheorem{remark}[theorem]{Remark}
\newtheorem{definition}[theorem]{Definition}

\newtheorem{construction}[theorem]{Construction}

\allowdisplaybreaks

\begin{document}
\title{An idempotent ring cannot be Morita equivalent to its ideal}
\address{$^1$Institute of Mathematics and Statistics, University of Tartu, Tartu, Estonia.}
\author{Kristo Väljako$^{1,2}$}
\address{$^2$Institute of Computer Science, University of Tartu, Tartu, Estonia.}
\date{\today}
\thanks{This research was supported by the Estonian Research Council grant PRG1204.}

\begin{abstract}
In this note it is proven that an idempotent ring cannot be Morita equivalent to its idempotent proper ideal.
\end{abstract}
\maketitle

\section{Preliminaries}

In this short note we study the Morita equivalence of idempotent rings. Note that we do not assume rings to have a unit element, unless specified otherwise. Recall that an associative ring $R$ is called \textbf{idempotent} if $RR = R$, i.e. for every $r \in R$ there exist $r_1,r_1',\ldots,r_k,r_k' \in R$ such that $r = r_1r_1' + \ldots + r_kr_k'$. We prove that if an idempotent ring $R$ is Morita equivalent to its idempotent ideal $S$, then $R = S$. This means that an idempotent ring cannot be Morita equivalent to its idempotent proper ideal. In the last section we consider two well-known examples of constructions for forming a unital ring $S$ from a ring $R$ such that the ring $R$ becomes an ideal of $S$.

Recall from \cite{Valjako}, the definition of an enlargement of a ring.

\begin{definition}[Def. 2.1 in \cite{Valjako}]
    We call a ring $R$ an \textbf{enlargement} of its subring $S$ if the conditions $R = RSR$ and
    $S = SRS$ hold. (We also say that $R$ is an \textbf{enlargement} of all rings isomorphic to such $S$.)
\end{definition}

Note that enlargements can be used to characterise the Morita equivalence of idempotent rings. Recall that a ring $T$ is called a \textbf{joint enlargement} of rings $R$ and $S$ if it is an enlargement of both $R$ and $S$.

\begin{theorem}[Theorem 3.9 in \cite{Valjako}]\label{TheoremMEiffEnlargements}
    Idempotent rings are Morita equivalent if and only if they have a joint enlargement.
\end{theorem}

\section{Result}

In this section we prove the main theorem of this short note. First we need two simple lemmas.

\begin{lemma}\label{LemIdealProperty}
    Let $R$ be ring and $S$ an idempotent ring. If $S \unlhd R$, then $RSR = S$.
\end{lemma}
\begin{proof}
    Let $R$ and $S$ be idempotent rings and $S \unlhd R$. By the definition of ideal, we have $RSR \subseteq S$. Due to $S$ being idempotent, we have $S = SS = SSS \subseteq RSR$. In conclusion, $S = RSR$.
\end{proof}

\begin{lemma}\label{LemEnlargementsOfIdealsEqual}
    Let $R$ and $S$ be idempotent rings and $S \unlhd R$. If $R$ and $S$ have a joint enlargement, then $R = S$.
\end{lemma}
\begin{proof}
    Let $T$ be the joint enlargement of $S$ and $T$. Then (upto isomorphism) we have $T = TST$ and $R = RTR$. Due to $S \unlhd R$ we also have $S = RSR$ (Lemma \ref{LemIdealProperty}). Now,
    \[
    R = RTR = R(TST)R = RT(RSR)TR = (RTR)S(RTR) = RSR = S.
    \]
    Hence we have $R = S$.
\end{proof}

Now we are ready to prove the main theorem of this note.

\begin{theorem}\label{TheoremMain}
    Let $R$ and $S$ be idempotent rings and $S \unlhd R$. If $R$ ans $S$ are Morita equivalent, then $R = S$.
\end{theorem}
\begin{proof}
    Let $R$ and $S$ be idempotent Morita equivalent rings, such that $S \unlhd R$. By Theorem \ref{TheoremMEiffEnlargements}, the rings $S$ and $R$ have a joint enlargement $T$. By Lemma \ref{LemEnlargementsOfIdealsEqual} we have $R = S$.
\end{proof}

\section{Corollaries}

In this section we recall two well-known constructions how to create a new unital ring from a ring $R$ such that $R$ will be an ideal of the new ring. By theorem \ref{TheoremMain}, neither of the rings obtained by these constructions can be Morita equivalent to the underlaying ring. We start with the Dorroh extension as introduced in \cite{Dorroh}.

\begin{construction}[Dorroh extension]
    Let $R$ be a ring. Consider the ring $R' = R \times \mathbb{Z}$, where addition is defined pointwise and multiplication as follows:
    \[
    (r,z)(s,w) := (rs + zs + wr,zw),
    \]
    where $(r,z),(s,w) \in R'$. It is easy to see that $R'$ is a unital ring with the unit element $(0,1)$.

    There exists a injective homomorphism of rings
    \[
    \iota\colon R \to R', \quad r \mapsto (r,0).
    \]
    One can easily check that $\Im\iota \cong R$ is an ideal of $R'$.
\end{construction}

\begin{corollary}
    An idempotent ring $R$ is not Morita equivalent to its Dorroh extension $R'$.
\end{corollary}

Next we recall the notion of a multiplier ring. There is a survey \cite{VanDeale} that collects many results concerning multiplier rings. Firs we say that a ring $R$ is called \textbf{non-degenarate} if for any $r \in R$ we have $r = 0$ if $rR = \{0\}$ or $Rr = \{0\}$.

\begin{construction}[Multiplier ring]
    Let $R$ be a non-degenerate ring. Consider the set
    \[
    \mathcal{M}(R) := \{(\rho,\lambda)\mid \rho,\lambda\colon R \to R,\ (\forall r,s \in R\colon\ \rho(r)s = r\lambda(s))\}.
    \]
    On the set $\mathcal{M}(R)$ define addition pointwise and multiplication as follows:
    \[
    (\rho,\lambda) * (\rho',\lambda') := (\rho' \circ \rho,\lambda \circ \lambda'),
    \]
    where $(\rho,\lambda),(\rho',\lambda') \in \mathcal{M}(R)$. It is easy to see that $\mathcal{M}(R)$ is a unital ring with the unit element $(\id,\id)$. It can be shown (Lemma 1.4 in \cite{VanDeale}) that if $(\rho,\lambda) \in \mathcal{M}(R)$, then $\rho\colon {_RR} \to {_RR}$ is a homomorphism of left $R$-modules and $\lambda\colon R_R \to R_R$ a homomorphism of right $R$-modules.

    Also, there exists a injective homomorphism of rings
    \[
    \iota\colon R \to \mathcal{M}(R), \quad r \mapsto (\rho_r,\lambda_r) = (\lbar r,r\lbar),
    \]
    where $\rho_r\colon r' \mapsto r'r$ and $\lambda_r\colon r' \mapsto rr'$. It can be shown that $\Im\iota \colon R$ is an ideal of $\mathcal{M}(R)$.
\end{construction}

\begin{corollary}
    An idempotent and non-degenrate ring $R$ in not Morita equivalent to its multiplier ring $\mathcal{M}(R)$.
\end{corollary}

Finally we make a little remark that connects Dorroh extension and multiplier rings.

\begin{remark}
    Let $R$ be a ring. Consider the category $\mathcal{C}$, which is a full subcategory of the category of unital rings $\mathsf{Ring}$ that consists of all rings to where $R$ can be inserted as an ideal. The category $\mathcal{C}$ has an initial and terminal object. The initial object of $\mathcal{C}$ is the Dorroh extension $R'$ and the terminal object is the multiplier ring $\mathcal{M}(R)$. Let $S$ be a unital ring such that $R \unlhd S$; then there exist unique homomorphisms of rings
    \begin{align*}
        &\mathrm{d}\colon\ R' \to S, \quad & (r,z) &\mapsto r + z1_S, \\
        &\mathrm{m}\colon\ S \to \mathcal{M}(R), \quad & s &\mapsto (\lbar s,s\lbar).
    \end{align*}
\end{remark}

\end{document}